\documentclass[scrbook]{ECOS_2023}
\usepackage{graphicx}
\usepackage{amsmath}
\usepackage{amsfonts}
\usepackage{amssymb}
\usepackage{mathtools}
\usepackage{enumitem}
\usepackage{verbatim}
\usepackage{hhline}
\usepackage{eurosym}
\usepackage{gensymb}
\makeatletter
\let\NAT@parse\undefined
\makeatother
\usepackage{hyperref}
\usepackage[title]{appendix}

\usepackage{booktabs}
\usepackage{multirow}

\usepackage[backend=bibtex,style=ieee,natbib=true]{biblatex}
\title{Reinforcement Learning for Joint Design and Control of Battery-PV Systems}

\author{Marine Cauz\textsuperscript{\,a,b}, Adrien Bolland\textsuperscript{\,c}, Bardhyl Miftari\textsuperscript{\,c}, Lionel Perret\textsuperscript{\,b}, Christophe Ballif\textsuperscript{\,a,d}, and Nicolas Wyrsch\textsuperscript{\,a}}

\address{
\textsuperscript{a} École polytechnique fédérale de Lausanne (EPFL), Institute of Electrical and Micro Engineering (IEM), Photovoltaics and Thin-Film Electronics Laboratory (PV-Lab), Neuchâtel, Switzerland. Marine.Cauz@epfl.ch, CA.
\and
\textsuperscript{b} Planair SA, Yverdon-les-bains, Switzerland
\and
\textsuperscript{c} University of Liège, Department of Electrical Engineering and Computer Science, Liège, Belgium
\and
\textsuperscript{d} Centre Suisse d'Electronique et de Microtechnique (CSEM), PV-Center, Neuchâtel, Switzerland}

\abstract{\normalsize
The decentralisation and unpredictability of new renewable energy sources require rethinking our energy system. Data-driven approaches, such as reinforcement learning (RL), have emerged as new control strategies for operating these systems, but they have not yet been applied to system design. This paper aims to bridge this gap by studying the use of an RL-based method for joint design and control of a real-world PV and battery system. The design problem is first formulated as a mixed-integer linear programming problem (MILP). The optimal MILP solution is then used to evaluate the performance of an RL agent trained in a surrogate environment designed for applying an existing data-driven algorithm. The main difference between the two models lies in their optimization approaches: while MILP finds a solution that minimizes the total costs for a one-year operation given the deterministic historical data, RL is a stochastic method that searches for an optimal strategy over one week of data on expectation over all weeks in the historical dataset. Both methods were applied on a toy example using one-week data and on a case study using one-year data. In both cases, models were found to converge to similar control solutions, but their investment decisions differed. Overall, these outcomes are an initial step illustrating benefits and challenges of using RL for the joint design and control of energy systems.

}


\begin{document}

\section{Introduction}
\label{sec:intro}

\subsection{Background and related work}

The current transition to renewable energy sources requires rethinking new energy systems, characterized by decentralized and intermittent production. The development of these systems typically occurs in two distinct steps, namely the design and control of these systems. The design problem involves identifying the design variables which are the optimal size of energy system components. The control problem aims to determine the control variables which are the optimal actions to operate the energy system components. Both design and control problem should jointly minimize a cost function and are typically solved sequentially. This paper explores the value of solving the design and control tasks, using a reinforcement learning (RL) method as appropriate design is intrinsically linked to subsequent operation. To evaluate the effectiveness of this approach, its performance are compared with that of the Mixed Integer Linear Programming (MILP) method. 

On the one hand, RL is a data-driven approach where an agent learns to make decisions in a dynamic environment through trial-and-error experience. It involves an agent interacting with an environment and receiving feedback in the form of rewards or penalties based on its actions, with the goal of maximizing its cumulative reward over time. On the other hand, Mixed Integer Linear Programming (MILP) is a mathematical optimization technique used to solve problems with linear constraints and integer variables. It involves formulating a mathematical model of the problem and using an optimization algorithm to find the best solution. Both RL and MILP methods will be used to benchmark the results of a one-year time series.

As highlighted in a recent review \citep{perera_applications_2021}, RL-based approaches have significant potential, yet not fully exploited, in the energy field. Specifically, the review points out that energy systems are typically designed using either MILP or heuristic methods, with RL approaches dedicated to their control. Integrating RL beyond energy flow control would open new interesting research questions. In \citep{perera_introducing_2020}, RL is used to support distributed energy system design due to its flexibility and model-free nature, which allows it to be adapted to different environments at different scales. However, they did not simultaneously address the dispatch and design problem as a distributed reward problem, as done in this work. Instead, they used a cooperative coevolution algorithm (COCE) to assist the optimization process. Jointly addressing the design and operation of energy systems is a key issue, especially for multi-energy systems, as discussed in \citep{fazlollahi_multi-objective_2013}, where multi-objective evolutionary algorithms (EMOO) and MILP are used to integrate biomass technologies in a multi-energy system. In \citep{majid_deep_2021}, the focus is on evolution algorithms and their comparison with deep reinforcement learning strategies. After clarifying the fundamental differences between the two approaches, the discussion revolves around their ability to parallelize computations, explore environments, and learn in dynamic settings. The potential of hybrid algorithms combining the two techniques is also investigated, along with their real-world applications.

RL-based frameworks are successfully applied to the operation of energy systems \citep{quest_3d_2022}, although these methods have not, to the authors' knowledge, been extended to solve real-world design problem in energy system. As reviewed in \citep{abdullah_reinforcement_2021}, RL-based frameworks are popular for addressing electric vehicle (EV) charging management, mostly with variants of the DQN algorithm, and outperform other traditional methods. In \citep{dorokhova_deep_2021}, various deep RL algorithms are benchmarked against rule-based control, model predictive control, and deterministic optimization in the presence of PV generation. The study, which aims to increase PV self-consumption and state-of-charge at departure, demonstrates the potential of RL for real-time implementation. For solving V2G control under price uncertainty, \citep{shi_real-time_2011} modeled the problem with a Markov Decision Process (MDP) \citep{uther_markov_2010}, a mathematical framework for modeling system where stochasticity is involved. Additionally, a linear MDP formulation is also used in \citep{manu_lahariya_computationally_nodate} to address the coordination of multiple charging points at once. Finally in \citep{sadeghianpourhamami_definition_2020}, a data-driven approach is defined and evaluated for coordinating the charging schedules of multiple EVs using batch reinforcement learning with a real use case. In conclusion, these studies provide valuable insights and tools for optimizing and improving energy systems, demonstrating the potential of RL to tackle the operation of complex energy systems.

\subsection{Contribution}

This work aims to evaluate the relevance of jointly designing and controlling an energy system using a deep RL approach. To achieve this purpose, two methods are benchmarked to address jointly the design and control problem of a real-world PV-battery system. The first method, MILP, computes the optimal design and control solution over a sequence of historical data. The second method, RL, computes the optimal design and a control policy through interactions with a simulator by trial and error. The specific RL algorithm used in this study is referred to as Direct Environment and Policy Search (DEPS) \citep{bolland_jointly_2022}. DEPS extends the REINFORCE algorithm \citep{williams_simple_1992} by combining policy gradient with model-based optimization techniques to parameterize the design variables. In this framework, an agent looks for the design and control variables that jointly maximize the expected sum of rewards collected over the time horizon of interest. The outcomes of both methods are discussed in the subsequent sections of this paper.

This paper is structured as follows. Section \ref{sec:methods} provides two formulations of the energy system, one designed for MILP and the other for RL, and discusses the methodology used to benchmark the results. In Section \ref{sec:results}, the outcomes of the study are presented, and these results are discussed in Section \ref{sec:discuss}, with a focus on the potential of RL for joint design and control of energy systems. Finally, the paper concludes with a summary in Section \ref{sec:ccl}.

\section{Method}
\label{sec:methods}
\subsection{Problem statement}
\begin{figure}[ht!]
    \centering
    \includegraphics[width=0.5\linewidth]{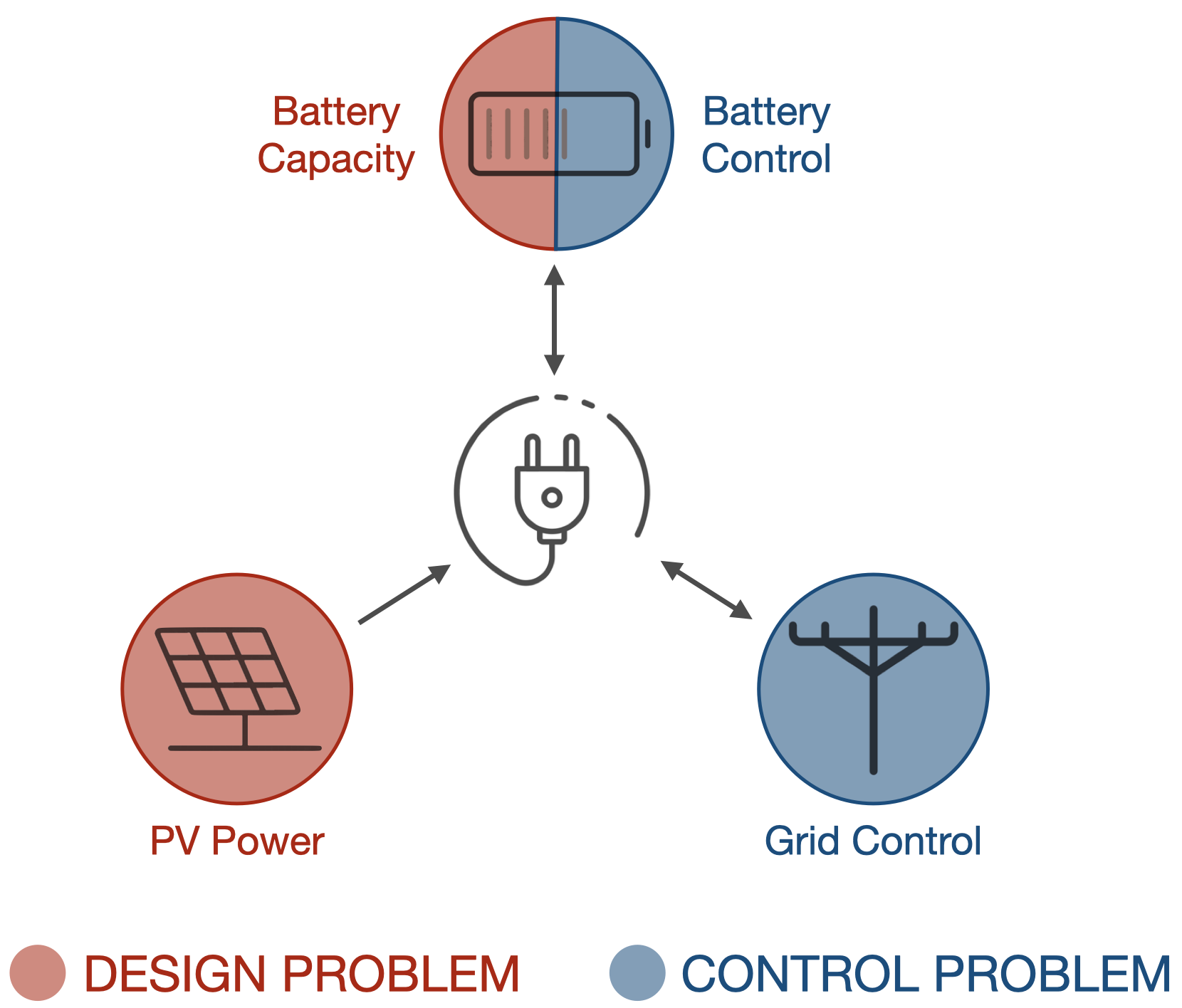}
    \caption{The energy system to be jointly designed and controlled is characterized by an electrical consumption, a battery, a photovoltaic system, and a grid connection. The design problem consists of determining the photovoltaic power and battery capacity, while the control problem aims to regulate the charging and discharging of the battery, as well as the import (resp. export) of electricity to (resp. from) the grid. The overall objective is to meet the electrical consumption needs while minimizing the costs of installing and operating the system.}
    \label{fig:system}
\end{figure}

The study is carried out for the energy system illustrated in Figure \ref{fig:system}, whose components are detailed in the subsections below. Overall, the system refers to an office building that has been fitted with a PV installation and a stationary lithium-ion battery to meet its own electricity consumption. Additionally, the building is connected to the electricity grid.

The objective of the study is to jointly propose a design of the PV and battery components, as well as a control strategy of the described energy system in order to minimize the total cost of its ownership. In the following Subsection \ref{sec:energy_sys}, the system is expressed as a mathematical program made-up of constraints and objectives. To be more precise, it is tackled as a Mixed-Integer Linear Program. Subsection \ref{sec:mdp} formulates a surrogate environment as an MDP. The latter represents the same dynamics and rewards as the original problem but the objective is to maximize the sum of rewards gathered over one week on expectation over the 52 weeks of the year of data. By doing so, it allows the use of the RL algorithm and expects the optimal solution to be close to the solution of the original problem. Results are discussed in Section \ref{sec:results}. Finally, for both methods, the energy system is studied over a finite time horizon $T$, on which all costs are evenly distributed across each time step $t$. The methodology and the context of the experiments conducted are specified in Subsection \ref{sec:benchmark}. \\

\subsection{Energy system}
\label{sec:energy_sys}
This subsection describes the physical constraints that apply to the components of the energy system. Theses components, in sequential order, consist of the PV panels, the battery, the electrical load and the power grid. The set of design and control variables and the parameters of the whole system, which is modeled as a discrete-time system, are gathered respectively in Table \ref{tab:variables} and \ref{tab:parameters}, respectively.

\begin{table}[ht!]
	\centering
	\small
	\begin{tabular}{@{}llllll@{}}
		\toprule
		& \textbf{Variable} & \textbf{Set} & \textbf{Unit} & \textbf{Description}  \\
		\midrule
		\parbox[t]{2mm}{\multirow{2}{*}{\rotatebox[origin=c]{90}{\textsc{Grid}}}}
		&$P^\textsc{imp}$ & $\mathbb{R}_+^T$ & kW & imported power (from the grid)\\
		&$P^\textsc{exp}$ & $\mathbb{R}_+^T$ & kW & exported power (to the grid)\\
		\midrule
		\parbox[t]{20mm}{\multirow{1}{*}{\rotatebox[origin=c]{90}{\textsc{PV}}}}
		&$P^\textsc{nom}$ & $\mathbb{R}_+$ & kW$_p$ & nominal power of the PV installation \\
		\midrule
		\parbox[t]{2mm}{\multirow{3}{*}{\rotatebox[origin=c]{90}{\textsc{Battery}}}}
		&$B$ & $\mathbb{R}_+$ & kWh & nominal capacity of the battery \\
		&$\textsc{soc}$ & $\mathbb{R}_+^T$ & kWh & state of charge of the battery \\
		&$P^B$ & $\mathbb{R}^T$ & kW & power exchanged with the battery \\
		\bottomrule
	\end{tabular}
    \caption{Set of design and control variables of the energy system studied.}
	\label{tab:variables}
\end{table}

\begin{table}[h!]
	\centering
	\small
		\begin{tabular}{@{}llllll@{}}
			\toprule
			& \textbf{Parameter} & \textbf{Value} & \textbf{Set} & \textbf{Unit}  & \textbf{Description}  \\
			\midrule
			\parbox[t]{2mm}{\multirow{4}{*}{\rotatebox[origin=c]{90}{\textsc{Grid}}}}
			&$C^\textsc{imp}_\textsc{grid}$ & 1 & $\mathbb{R}$ & \euro{}/kWh & imported electricity price \\ 
			&$C^\textsc{exp}_\textsc{grid}$ & -0.05 & $\mathbb{R}$ & \euro{}/kWh & exported electricity price \\
		    &$C_\textsc{grid}$ &  & $\mathbb{R}^T$ & \euro{} & electricity network cost \\
		    &$P_\textsc{grid}^{\textsc{max}}$ & 10'000 & $\mathbb{R}_+$ & kW & grid connection power \\
			\midrule
			\parbox[t]{2mm}{\multirow{13}{*}{\rotatebox[origin=c]{90}{\textsc{PV}}}}
			&$P^\textsc{nom}_\textsc{min}$ & 0 & $\mathbb{R}_+$ & kW$_p$  & minimal nominal PV power \\ 
			&$P^\textsc{nom}_\textsc{max}$ & 200 & $\mathbb{R}_+$ & kW$_p$  & maximal nominal PV power \\			
			&$P^\textsc{prod}$ &  & $\mathbb{R}_+^T$ & kW & generated PV power \\
    		&$\overline{P}^\textsc{prod}$ &  & $\mathbb{R}_+^T$ & kW & expected generated PV power \\
            &$\overline{p}^\textsc{prod}$ &  & $\mathbb{R}_+^T$ & kW & normalised PV power \\
    	  &$L^\textsc{pv}$  & 20 & $\mathbb{N}$ & years & PV lifetime\\
            &$\textsc{R}_\textsc{pv}$ &  & $\mathbb{R}_+$ & - & annuity factor\\ 
			&$\textsc{ox}_\textsc{pv}^\textsc{fix}$  & 3 & $\mathbb{R}_+$ & \euro{} & \textsc{opex} PV fixed cost\\  
			&$\textsc{ox}_\textsc{pv}^\textsc{var}$ & 10 & $\mathbb{R}_+$ & \euro{}/kW & \textsc{opex} PV variable cost\\
			&$\textsc{cx}_\textsc{pv}^\textsc{fix}$ & 50 & $\mathbb{R}_+$ & \euro{} & \textsc{capex} PV fixed cost\\
			&$\textsc{cx}_\textsc{pv}^\textsc{var}$ & 200 & $\mathbb{R}_+$ & \euro{}/kW & \textsc{capex} PV variable cost\\ 
			\midrule
			\parbox[t]{2mm}{\multirow{9}{*}{\rotatebox[origin=c]{90}{\textsc{Battery}}}}
            &$B_\textsc{min}$ & 0 & $\mathbb{R}_+$ & kWh  & minimal nominal battery capacity \\ 
			&$B_\textsc{max}$ & 200 & $\mathbb{R}_+$ & kWh  & maximal nominal battery capacity \\
			&$\eta^\textsc{b}$ & 0.9 & $\left]0,1\right]$ & - & battery  efficiency\\
            &$L^\textsc{b}$ & 30 & $\mathbb{N}$ & years & battery lifetime\\
            &$\textsc{R}_\textsc{B}$ &  & $\mathbb{R}_+$ & - & annuity factor\\
            &$\textsc{ox}_\textsc{b}^\textsc{fix}$ & 5 & $\mathbb{R}_+$ & \euro{} & \textsc{opex} Battery fixed cost\\  
			&$\textsc{ox}_\textsc{b}^\textsc{var}$ & 6 & $\mathbb{R}_+$ & \euro{}/kW & \textsc{opex} Battery variable cost\\
			&$\textsc{cx}_\textsc{b}^\textsc{fix}$ & 30 & $\mathbb{R}_+$ & \euro{} & \textsc{capex} Battery fixed cost\\
			&$\textsc{cx}_\textsc{b}^\textsc{var}$ & 110 & $\mathbb{R}_+$ & \euro{}/kW & \textsc{capex} Battery variable cost\\ 
			\midrule
			\parbox[t]{2mm}{\multirow{6}{*}{\rotatebox[origin=c]{90}{\textsc{System}}}} 
			&$T$  &  & $\mathbb{N}$& -& time horizon \\ 
			&$\Delta t$ & 1 &  $\mathbb{R}_+^T$ & h  & time steps \\
			&$h_t$ &  & $[0:23]$ & h & hour of the time step \\
			& $r$ & 0.05 & $\mathbb{R}$ & -  & discount rate \\ 
			& $P^\textsc{load}$ &  &  $\mathbb{R}_+^T$ & kW & uncontrollable electricity consumption\\
			& $\overline{P}^\textsc{load}$ &  & $\mathbb{R}_+^T$ & kW & expected electricity consumption\\
			\bottomrule
		\end{tabular}
    \caption{Set of parameters of the energy system studied.}
	\label{tab:parameters}
\end{table}

\subsubsection*{PV system}
The objective of the PV installation is to generate electricity on-site to fulfill the local electricity demand. The design of this component is one of the two design variables that will result from the optimization process. The range of the suitable nominal power $P^\textsc{nom}$, corresponding to its design variable, is set in Eq. \eqref{eq:pv_nom} and the production at time $t$ is directly proportional to this nominal design variable as shown in Eq. \eqref{eq:pv_prod}. The normalized annual curve $\overline{p}_t^\textsc{prod}$ corresponds to the actual hourly averaged PV production power of the building. 

    \begin{align}
        \label{eq:pv_nom}
        P^\textsc{nom}_\textsc{min} &\le P^\textsc{nom} \le P^\textsc{nom}_\textsc{max} \\
        \label{eq:pv_prod}
            P^\textsc{prod}_t & = P^\textsc{nom} \cdot \overline{p}^\textsc{prod}_t 
        \end{align}
        
The \textsc{capex} and \textsc{opex} values, which are respectively the initial investment and the annual maintenance cost, of the installation are made up of a fixed and a variable part to take account of potential scale effects.
    \begin{align}
        \label{eq:pv_cx}
        \textsc{cx}_{\textsc{pv}} &= \textsc{cx}_{\textsc{pv}}^\textsc{fix} + 
        \textsc{cx}_{\textsc{pv}}^\textsc{var} \cdot P^\textsc{nom}\\
        \label{eq:pv_ox}
        \textsc{ox}_{\textsc{pv}} &= \textsc{ox}_{\textsc{pv}}^\textsc{fix} + 
        \textsc{ox}_{\textsc{pv}}^\textsc{var} \cdot P^\textsc{nom}
    \end{align}

\subsubsection*{Battery}
To maximize the potential for on-site self-consumption, a stationary lithium-ion battery is available. The design of this component, corresponding to its capacity $B$, is the second design variable to determine during the optimization process. This battery capacity can vary in the range of Eq. \eqref{eq:bat_capacity}.

    \begin{align}
        B_\textsc{min} &\le B \le B_\textsc{max}
        \label{eq:bat_capacity}
    \end{align}
    
The state of charge variable, $\textsc{soc}_t$, changes as a function of the power exchanged with the battery denoted $P^B_t$. This power is constrained, for charging, by the nominal capacity, Eq. \eqref{eq:bat_ch}, and, for discharging, by the energy stored, Eq. \eqref{eq:bat_disc}. Additionally, the battery efficiency, denoted $\eta^\textsc{b}$, is assumed identical for both the charging and the discharging processes. 

    \begin{align}
        \label{eq:bat_ch}
        P^B_t &\le \frac{\textsc{B} - \textsc{soc}_t}{\Delta t} &\text{if } P^B_t \ge 0\\
        P^B_t &\ge \frac{-\textsc{soc}_t}{\Delta t} &\text{if } P^B_t \le 0
        \label{eq:bat_disc}
    \end{align}
    
Knowing the power exchanged with the battery, the state of charge can be updated:
    \begin{align}
        \textsc{soc}_{t+1} &= 
        \begin{cases}
            \textsc{soc}_t + P^B_t \cdot \eta^\textsc{b} \cdot \Delta t &\text{ if } P^B_t \ge 0 \\
            \textsc{soc}_t + \frac{P^B_t}{\eta^\textsc{b}} \cdot \Delta t &\text{ if } P^B_t < 0
        \end{cases}
        \label{eq:bat_soc}
    \end{align}
At the beginning of the optimization, i.e., $t=0$, the battery \textsc{soc} is set to half of its capacity value, to initialize the model. Moreover, to avoid any artificial benefit, the final \textsc{soc} is constrained to be equal to the initial value, as formulated in Eq. \eqref{eq:soc_init}.
    \begin{align}
        \textsc{soc}_{t=0} &= \frac{B}{2}\\
        \textsc{soc}_{t=0} &= \textsc{soc}_{t=T}
    \label{eq:soc_init}
    \end{align}
  
Similar to the PV plant, the \textsc{capex} and \textsc{opex} of the battery consist of both fixed and variable parts.
    \begin{align}
        \label{eq:bat_cx}
        \textsc{cx}_{\textsc{B}} &= \textsc{cx}_{\textsc{B}}^\textsc{fix} + 
        \textsc{cx}_{\textsc{B}}^\textsc{var} \cdot \textsc{B}\\
        \label{eq:bat_ox}
        \textsc{ox}_{\textsc{B}} &= \textsc{ox}_{\textsc{B}}^\textsc{fix} + 
        \textsc{ox}_{\textsc{B}}^\textsc{var} \cdot \textsc{B}
    \end{align}

\subsubsection*{Electrical load}
The electrical load used in this project is real data from an office building in Switzerland. This consumption is monitored on an hourly basis and reflects the consumption patterns of office days. This building load power, $P^\textsc{load}_t$, is provided as input and corresponds to an actual measurement sampled by hours over a year.

\subsubsection*{Electrical grid}
To absorb excess solar production or to meet the electricity consumption in the absence of local production, the system is connected to the low-voltage electrical grid. This connection is modeled here as a single balance equation, called the conservation of electrical power, shown in Eq. \eqref{eq:energy_balance}. The power imported from the grid is referred to as $P^\textsc{imp}_t$ and the power injected is referred to as $P^\textsc{exp}_t$.
    \begin{align}
        P^\textsc{prod}_t + P^\textsc{imp}_t = P^\textsc{load}_t + P^\textsc{B}_{t} + P^\textsc{exp}_t
    \label{eq:energy_balance}
    \end{align}
The  grid power value at each time $t$ is derived from Eq. \ref{eq:energy_balance}, and the power limit can be described as follows.
    \begin{align}
         0 &\le P^\textsc{imp}_t \le P^\textsc{max}_\textsc{grid} \\
         0 &\le P^\textsc{exp}_t \le P^\textsc{max}_\textsc{grid}
    \end{align}
Based on the import and export power, the total cost of supplying electricity through the network $C_\textsc{grid}$ can be computed.
    \begin{equation}
            C_\textsc{grid} = \sum_{t=0}^{T-1}  C_{\textsc{grid},t} = \sum_{t=0}^{T-1} P^\textsc{imp}_t \cdot C^\textsc{imp}_{\textsc{grid},t} - P^\textsc{exp}_t \cdot C^\textsc{exp}_{\textsc{grid},t}
            \label{eq:grid_cost}
    \end{equation}

\subsubsection*{Objective function}
The objective of this study is to propose a design for the PV and battery components, along with their dispatch, with the aim of minimizing the total cost of ownership. This objective function, of minimizing the overall cost of the system, can be formulated as follows.
    \begin{equation}
    \min \ \ \textsc{totex}
    \label{eq:milp_obj}
    \end{equation}
    
The total cost of the system, denoted \textsc{totex}, is composed of the \textsc{capex} and \textsc{opex} of both PV and battery components, as well as the grid cost. 
    \begin{align}
    \label{eq:totex}
    \textsc{totex} &= \textsc{opex} + \textsc{capex} + C_\textsc{grid} \\
    \label{eq:OPEX}
    \textsc{opex}& = \textsc{ox}_\text{pv} + \textsc{ox}_\textsc{B} \\
    \label{eq:CAPEX}
    \textsc{capex} &= \textsc{cx}_\text{pv} \cdot R_{pv} + \textsc{cx}_\textsc{B} \cdot R_{B}
    \end{align}
The \textsc{opex} and grid cost are computed over a finite time period $T$. However, the \textsc{capex} is an investment cost that is independent of $T$. To enable the adaptation of the investment cost to the project duration, an annuity factor $R$ adjusts the \textsc{capex} for the finite time horizon $T$. This annuity factor is computed according to Eq. \eqref{eq:Rfact}, by taking into account the values of $T$, the annual discount rate $r$, and the lifetime $L$ of the component. This formula includes a scaling factor $\frac{T}{8760}$ to adapt $R$ to the period $T$, based on the assumption that $T$ is expressed in hours since 8760 is the number of hours in a year. 
    \begin{align}
    \label{eq:Rfact}
    R& = \frac{r \cdot (1 + r) ^L}{ (1 + r) ^ L - 1} \cdot \frac{T}{8760}
    \end{align}

\subsection{MDP formulation}
\label{sec:mdp}
This section presents an alternative formulation of the problem as a Markov Decision Process (MDP), which is a well-established framework for modeling sequential decision-making problems. This alternative formulation is required for applying DEPS. More precisely, an MDP$(S, A, P, R, T)$, as presented below, consists of the following elements: a finite set of states $S$, a finite set of actions $A$, a transition function $P$, a rewards function $R$, and a finite time horizon $T$.

\subsubsection*{State Space}
The state of the system can be fully described by
    \begin{align}
        s_t &= (h_t, d_t, \textsc{soc}_t,  \overline{P}^\textsc{prod}_t,  \overline{P}^\textsc{load}_t)\\
        &\in S = \{0, ..., 23\} \times \{0, ..., 364\} \times [0, B] \times \mathbb{R}_+ \times \mathbb{R}_+ 
    \end{align}
    \begin{itemize}
        \item $h_t \in \{0, ..., 23\}$ denotes the hour of the day at time $t$. The initial value is set to 0. 
        \item $d_t \in \{0, ..., 364\}$ denotes the day of the year at time $t$. The initial value is set randomly.
        \item $\textsc{soc}_t$ is the state of charge of the battery at time $t$, this value is upper bounded by the nominal capacity of the installed battery $\textsc{B}$. The value is set initially to a random value during the training process and to half of its capacity during the validation process.
        \item $\overline{P}^\textsc{prod}_t$ represents the expected PV power at time $t$. This value is obtained by scaling normalized historical data $\overline{p}^\textsc{prod}_t$ with the total installed PV power ($P^\textsc{nom}$) and considering $h_t$ and $d_t$ values.
        \item $\overline{P}^\textsc{load}_t$ denotes the expected value of the electrical load at time $t$. The load profile is determined using historical data that corresponds to the same hour and day as the PV power.
    \end{itemize}

\subsubsection*{Action Space}
The action of the system corresponds to the power exchanged with the battery.
\begin{align}
        \widetilde{a_{t}} &= (\widetilde{P}^B_t)
\end{align}
After projecting the action to fall within the acceptable range specified by Eq. \eqref{eq:bat_ch} and \eqref{eq:bat_disc}, the resulting value is used as $a_{t}$, as shown in Eq. \eqref{eq:act}. This corresponds to the power exchanged with the battery, denoted $P^B_t$, this value is positive when the battery is being charged and negative when it is being discharged.
    \begin{align}
        P^B_t = 
        \begin{cases}
          \frac{\textsc{B} - \textsc{soc}_t}{\Delta t} &\text{ if } \widetilde{P}^B_t > \frac{\textsc{B} - \textsc{soc}_t}{\Delta t}\\
          \frac{\textsc{soc}_t}{\Delta t} &\text{ if } \widetilde{P}^B_t < - \frac{\textsc{soc}_t}{\Delta t} \\
          P^B_t &\text{otherwise}
        \end{cases}
    \label{eq:act}
    \end{align}

\subsubsection*{Transition Function}
Each time step $t$ in the system corresponds to one hour, which implies the evolution specified in Eq. \eqref{eq:h+1} of the state variable $h$ and every 24 time steps, the day is incremented by 1.
    \begin{align}
        \label{eq:h+1}
        h_{t+1} &= (h_t + 1) \text{ mod } 24\\
        d_{t+1} &= \text{Int}(\frac{h_t + 1}{24})
        \label{eq:d+1}
    \end{align}
    \noindent where the function $Int$ takes the integer value of the expression in Eq. \eqref{eq:d+1}.

The $\textsc{soc}_t$ of the battery is updated as Eq. \eqref{eq:bat_soc}, based on the projected action value, and all other state variables are taken from input data.
    \begin{align}
        \overline{P}^\textsc{prod}_{t+1} &= \overline{p}^\textsc{prod}_{h_{t+1}, d_{t+1}} \cdot P^\textsc{nom}\\
        \overline{P}^\textsc{load}_{t+1} &= \overline{p}^\textsc{load}_{h_{t+1}, d_{t+1}}
    \end{align}
 
\subsubsection*{Reward Function}
The reward signal to optimize the agent's actions in RL serves a similar aim as the objective function in the MILP formulation. Therefore, the reward here is the opposite value of the \textsc{TOTEX} defined at Eq. \eqref{eq:totex}. This cost is composed of (i) the investment cost, (ii) the operating cost and (iii) the cost from the purchase and resale of electricity from the grid defined in Eq. \eqref{eq:grid_cost}.
    \begin{align}
        r_t &= - \textsc{totex}_t \\
        &= - \textsc{capex} - \textsc{opex} - C_{\textsc{grid},t} \\
        &= - \textsc{capex} - \textsc{opex} - P^\textsc{imp}_t \cdot C^\textsc{imp}_{\textsc{grid},t} + P^\textsc{exp}_t \cdot C^\textsc{exp}_{\textsc{grid},t}
        \label{eq:reward}
    \end{align}
    \noindent where the grid cost is the only time-dependent factor, while \textsc{capex} and \textsc{opex} are fixed values for a specific value of $P^{\textsc{nom}}$ and $B$.\\

\subsection{Methodology}
\label{sec:benchmark}

This subsection discusses the fundamental differences between the two methods (i.e., MILP and RL), along with the experimental protocol that was employed to compare the results. As discussed briefly earlier, although the same problem is aimed to be solved, the methods under study are fundamentally different.

MILP is a method for solving problems that involves optimizing a linear function of variables that are either integer or constrained by linear equalities, as the problem described in Subsection \ref{sec:energy_sys} The MILP algorithm solves the optimization problem by iteratively adjusting the values of the design and control variables, subject to the constraints, until it finds the optimal solution that maximizes or minimizes the objective function, depending on the problem's goal. This method is applied to the problem described in Subsection \ref{sec:energy_sys} over a one-year time horizon ($T=8760$). The solution is said to be computed with perfect foresight meaning that all variables are selected accounting for the future realization of (normally unknown) events in the time series, providing an optimistic upper bound on the true performance of the control and design. Concretely, the MILP problem is here encoded in the Graph-Based Optimization Modelling Language (GBOML) \citep{miftari_gboml_2022} paired with the Gurobi solver \citep{gurobi_gurobi_2020}.

In contrast, RL is a stochastic optimization method that learns from experience through trials and errors. In this study, we use DEPS \citep{bolland_jointly_2022}, an algorithm optimizing design and control variables in an MDP, as the one described in Subsection \ref{sec:mdp}, with a finite-time horizon. The agent receives feedback in the form of rewards when it selects a particular design and performs specific actions. The objective of the agent is to maximize the expected cumulative reward, which drives it to learn a design and a control policy.  Ideally, as with MILP, the time horizon should be annual, or cover the entire lifetime of the system, taking into account seasonal production and consumption fluctuations and/or equipment aging. However, such extended time horizons are unsuitable for this RL approach. Therefore, to strike a balance between a horizon that is short enough for DEPS and long enough to observe the consequences of decisions on the system, a horizon of $T=168$ hours, i.e., 7 days, is defined. Additionally, for each simulation, the initial day is sampled uniformly from the year-long data set and the initial state-of-charge of the battery is also sampled uniformly at random. As the reward is optimized on expectation over all days, the resulting design and control policy is expected to account for the seasonality and other different hazard in the historical data. The DEPS algorithm is trained on a predetermined number of iterations. The PV power and battery capacity values obtained from the last iteration of the algorithm are then taken as the values of the design variables and the final policy is used for the control.

Unlike MILP, the RL method does not secure optimality, therefore the experimental protocol aims to compare both results to see how far the RL solution is from the optimal one.  The experimental protocol is conducted in two distinct scenarios to differentiate the impact of adding the design variables in the joint problem. The first control-only scenario (\textsc{CTR}) assessed the control variables when the design variables are fixed. The second scenario, considering both the control and design (\textsc{CTR} \& \textsc{DES}) problem, allows for flexibility in designing the battery capacity and PV power, the two design variables. To benchmark the performance of both methods in each scenario (i.e., \textsc{CTR} and \textsc{CTR} \& \textsc{DES}), the reward and income value are reported. The reward value is computing according to Eq. \eqref{eq:reward} for the RL method. To estimate the average reward value for the MILP method, all reward values $r_t$ are averaged over time horizons of $T=168$. Comparing the average cumulative reward value of the MILP method to that of the RL method  provides a first benchmark for evaluating the performance of both approaches. However, as shown in Eq. \eqref{eq:reward}, only the grid cost is time-dependent, while the \textsc{capex} and \textsc{opex} depend solely on investment decisions. Therefore, the income value is defined as the average reward value, but it only includes the grid cost and can be computed as follows:
    \begin{align}
        Income &= \sum_{t=0}^{T-1} - P^\textsc{imp}_t \cdot C^\textsc{imp}_{\textsc{grid},t} + P^\textsc{exp}_t \cdot C^\textsc{exp}_{\textsc{grid},t}
    \end{align}

Finally, the experiments are performed in two steps. First, to perform a simple comparative study, working on a same finite time horizon $T=168$, both methods are conducted using data from a single summer week. Second, the data set is extended to include the one-year data set.

\section{Results}
\label{sec:results}

The energy system presented in Section \ref{sec:methods} is solved using the RL and MILP approaches with parameter values listed in Table \ref{tab:parameters}. To differentiate the performance of the DEPS algorithm for control and design aspects, the study is conducted in two distinct scenarios. The first control-only scenario (\textsc{CTR}) assessed the control aspect for fixed design variables, meaning that the PV power and battery capacity are fixed. The second scenario, considering both the control and design (\textsc{CTR} \& \textsc{DES}) aspects, allows for flexibility in designing the battery capacity and PV power. The two following subsections describe the results of the study performed in two steps, over the one-week and one-year data set, respectively.\\

\subsection{A one-week toy example}

In order to perform a simple comparative study, both \textsc{CTR} and \textsc{CTR} \& \textsc{DES} analyses were conducted using data from a single summer week. This enables to optimize both methods on the same time horizon. This means training the RL algorithm on the same $168$ time steps, with an initial day uniformly randomly selected over the week but an initial hour fixed at midnight. Additionally, during the training phase, the battery's initial \textsc{soc} is uniformly sampled such that the RL algorithm is presented with a large variety of scenarios for improving the quality of the learned policy. The results for both the \textsc{CTR} scenario, where the design variables (i.e., the PV power and battery capacity) are fixed, and the \textsc{CTR} \& \textsc{DES} scenario, where the PV and battery design variables are optimized in addition to control, are presented in Table \ref{tab:results_week}.

\begin{table}[ht!]
	\centering
	\small
	\begin{tabular}{cccccc}
		\toprule
		& &\textbf{Unit} & \textbf{Optimal RL} & \textbf{Optimal MILP}& \textbf{MILP solution based}  \\
  		&  & & \textbf{solution} & \textbf{solution}& \textbf{on RL design}  \\
		\midrule
		\parbox[t]{2mm}{\multirow{3}{*}{\rotatebox[origin=c]{90}{\textsc{CTR}}}}
		&T &hours &168 &168 &\\
 &Reward &\euro{} &-66 &-66 &\\
   &Income &\euro{} &-30 &-30 & \\
		\midrule
		\parbox[t]{2mm}{\multirow{5}{*}{\rotatebox[origin=c]{90}{\textsc{CTR \& DES}}}}
		&T &hours &168 &168 &168 \\
  &Reward &\euro{} &-40 &-46 &-53 \\
		&Income &\euro{} &-4 &0 &-17 \\
		&Battery capacity &kWh &62 &40 &62 \\
		&PV power &kWp &41 &103 &41 \\
		\bottomrule
	\end{tabular}
    \caption{Results of RL and MILP solutions on one-week data for control-only (\textsc{CTR}) and control and design (\textsc{CTR} \& \textsc{DES}) scenarios. \textit{T} denotes the time horizon in hours, while \textit{income} represents the cost of buying and selling electricity from the grid, \textit{reward} is the average cumulative reward value, and \textit{battery capacity} and \textit{PV power} indicate the values of design variables, which were set to 31.89 kWh and 55.81 kWp, respectively, in the \textsc{CTR} scenario. In the RL model, reward and income values were obtained by reloading the trained model with the determined design variables. Results were computed using an initial state of charge of the battery set to 50\% of its capacity for both models. However, the RL model does not take into account the Eq. \eqref{eq:soc_init}.}
    \label{tab:results_week}
\end{table}

\subsubsection{RL and MILP optimal objective values are similar in both scenarios but with different designs in the control and design scenario.}

Table \ref{tab:results_week} shows that in the \textsc{CTR} scenario, the results of the RL approach are similar to those of MILP. This confirms that the DEPS algorithm is able to converge to the optimal solution of this specific problem. In the \textsc{CTR} \& \textsc{DES} scenario, RL design variables differ from the MILP solution, resulting in an unexpected higher reward value (-40) than the MILP optimal one (-46). A detailed analysis reveals that this unexpected higher value is due to Eq. \eqref{eq:soc_init}, which is not imposed in the MDP. In order to validate this analysis, the additional grid cost needed to fulfill Eq. \eqref{eq:soc_init} has been computed, taking into account the battery's final \textsc{soc} obtained with the RL approach. As a result, the reward value has increased to -67 (instead of -40). This clearly highlights the importance that Eq. \eqref{eq:soc_init} plays in term of overall objective.

\subsubsection{The \textsc{CTR} \& \textsc{DES} scenario highlights differences in RL and MILP strategies.}

It is seen from Table \ref{tab:results_week} that in the second scenario, the optimal design variables of the RL and MILP solutions differ. Finding different values in design variables shows that the DEPS algorithm is able to identify solutions with comparable reward but using different design strategies. In order to study the sensitivity of the optimal solution, the MILP method was applied by imposing the design variable values obtained with the RL, as it can be seen in the last column of Table \ref{tab:results_week}. This indicates that the RL design solution is less optimal (-53) than the MILP one (-46).

\subsection{A one-year case study}
Optimal solutions of RL and MILP methods in both scenarios are now computed using data from a full year. The time horizon for the RL algorithm is still equal to $T=168$, but the starting days are uniformly randomly selected over the year. The RL algorithm is trained over a pre-determined number of 100'000 iterations and the values of the RL design variables considered are the ones from the final iteration.  The results are shown in Table \ref{tab:results_year}.

\begin{table}[h!]
	\centering
	\small
	\begin{tabular}{@{}cccccc@{}}
		\toprule
		&  &\textbf{Unit} & \textbf{Optimal RL } & \textbf{Optimal MILP } & \textbf{MILP solution based}  \\
  &  & & \textbf{solution} & \textbf{solution} & \textbf{on RL design}  \\
		\midrule
		\parbox[t]{2mm}{\multirow{3}{*}{\rotatebox[origin=c]{90}{\textsc{CTR}}}}
		&T &hours &168 &8760 &\\
		&Reward &\euro{} &-268 &-228 &\\
  	&Income &\euro{} &-220 &-196 & \\
		\midrule
		\parbox[t]{2mm}{\multirow{5}{*}{\rotatebox[origin=c]{90}{\textsc{CTR \& DES}}}}
		&T &hours &168 &8760 &8760 \\
		&Reward &\euro{} &-250 &-205 &-247 \\
		&Income &\euro{} &-208 &-164 &-218  \\
		&Battery capacity &kWh &44 &95 &44 \\
		&PV power &kWp &57 &81 &57 \\
		\bottomrule
	\end{tabular}
    \caption{Results of RL and MILP solutions on one-year data for \textsc{CTR} and \textsc{CTR} \& \textsc{DES} scenarios. \textit{T} denotes the time horizon in hours, while \textit{income} represents the cost of buying and selling electricity from the grid, \textit{reward} is the (expected) cumulative reward value, and \textit{battery capacity} and \textit{PV power} indicate the design variable values, which were set to 64.9 kWh and 63.65 kWp, respectively, in the \textsc{CTR} scenario. The RL solution is based on the trained model to determine the reward and income values, based on an average of 1'000 simulations over $T=168$. The MILP solution is computed for a one-year time horizon ($T=8760$). Both models use an initial state of charge (\textsc{soc}) of the battery set to 50\% of its capacity. However, the MILP model has an additional constraint specifies in Eq. \eqref{eq:soc_init}.}
    \label{tab:results_year}
\end{table}

\subsubsection{The difficulty of generalizing a policy with stochasticity in the model and on the estimation of the expectation}

\begin{figure}[h!]
    \centering
    \includegraphics[width=1\linewidth]{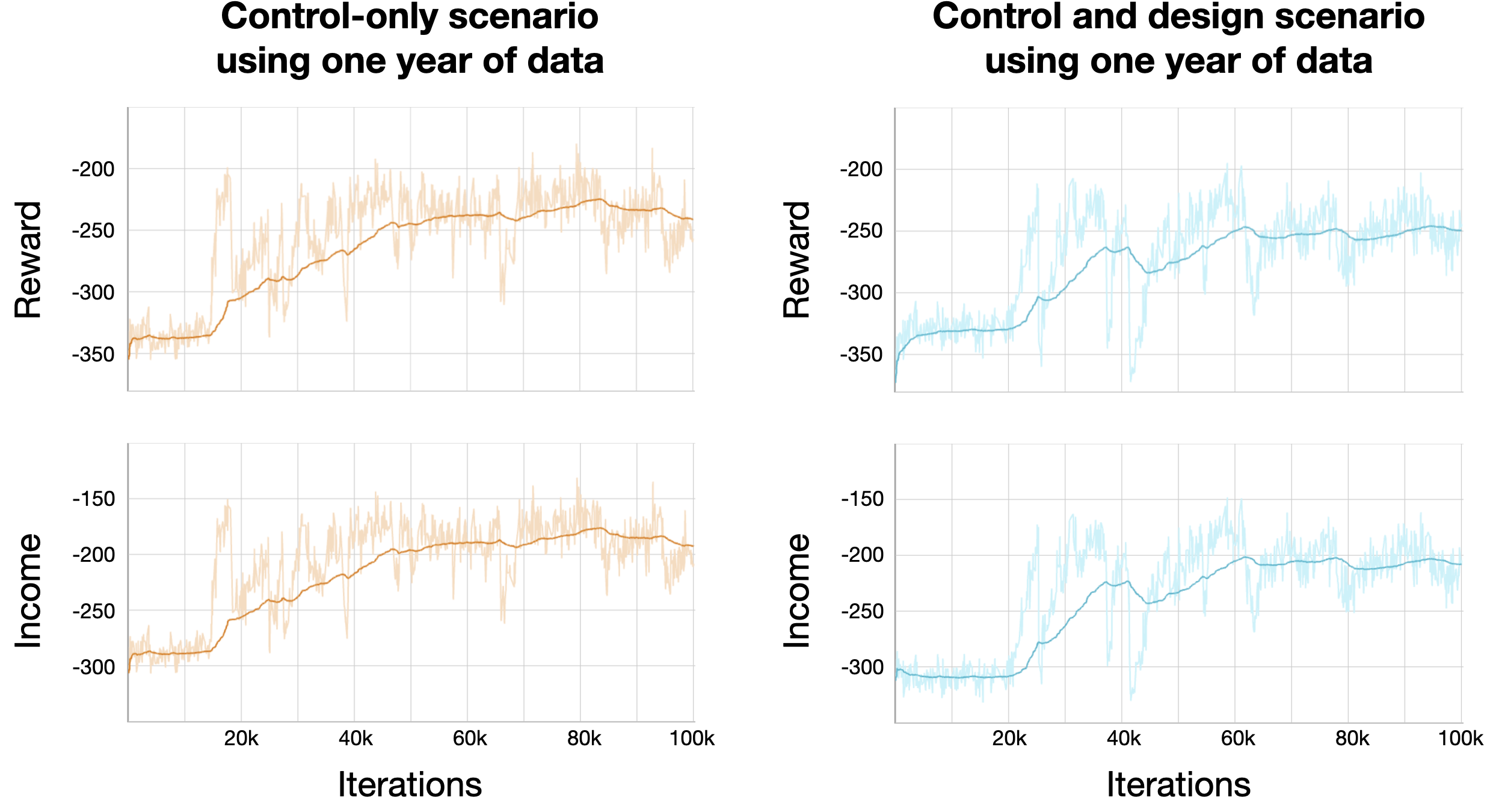}
    \caption{Value of reward and income obtained by the DEPS algorithm at each iteration for both scenarios. The left plots show the reward and income values for the \textsc{CTR} scenario, while the right plot shows the same values for the \textsc{CTR} \& \textsc{DES} scenario. The light curve shows the exact values for each time step, while the dark curve displays the corresponding smoothed values. In the \textsc{CTR} scenario, the difference between the reward and income values remains constant at 31.93 due to fixed design variables, with a battery size of 64.9 kWh and a PV power of 63.65 kWp. However, in the \textsc{CTR} \& \textsc{DES} scenario, the battery size and PV power output vary from 0 to 200 kWh and kWp, respectively.}
    \label{fig:ctr_year}
\end{figure}

It can be seen from Table \ref{tab:results_year} that in both the \textsc{CTR} and \textsc{CTR} \& \textsc{DES} scenarios, the optimal reward obtained by the RL method is poorer than the MILP optimal rewards. Furthermore, as depicted in Fig. \ref{fig:ctr_year}, due to the significant variations in the input data, the reward and income values exhibit substantial fluctuations across iterations.

During training in the \textsc{CTR} scenario (Fig. \ref{fig:ctr_year}, left), the RL model achieved maximum reward and income values of -180 and -131, respectively, which are significantly better than the final results obtained from both methods in Table \ref{tab:results_year}. This could suggest that depending on the set of weeks that are averaged at each iteration, it is possible to obtain a better or worse reward. Therefore, it seems important to work with a sufficiently representative number of weeks throughout the year. A similar observation can be made in the \textsc{CTR} \& \textsc{DES} scenario, where the maximum reward and income values achieved were -195 and -148, respectively (Fig. \ref{fig:ctr_year}, right).

\subsubsection{The RL method seems to promote lower design variable values}

\begin{figure}[h!]
    \centering
    \includegraphics[width=1\linewidth]{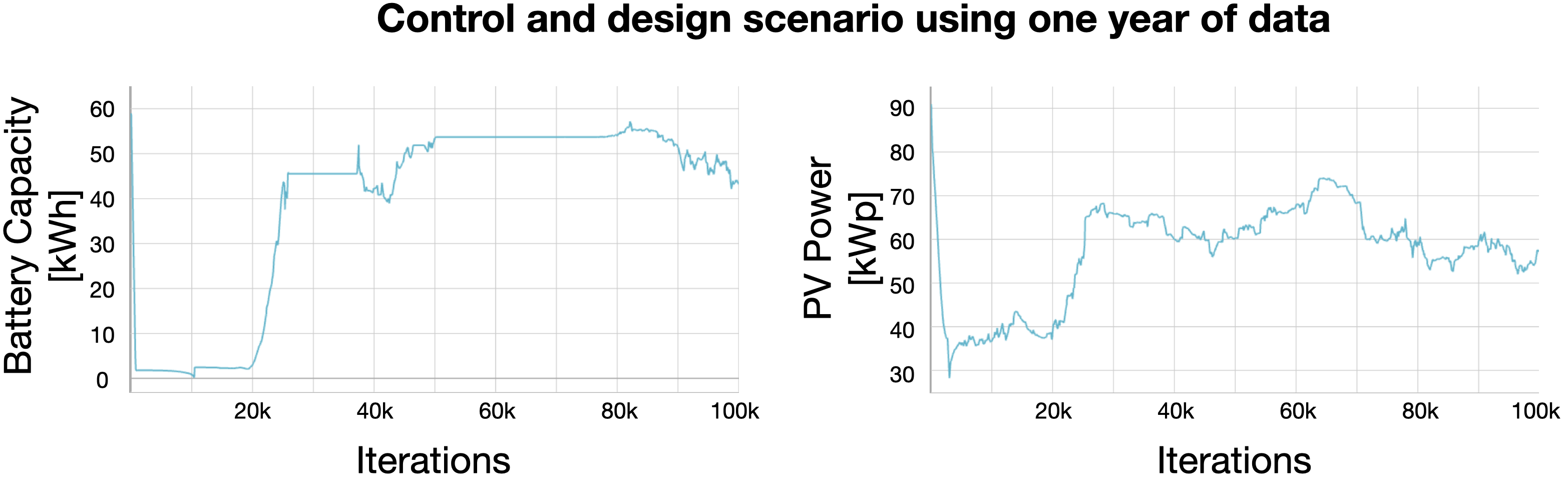}
    \caption{Value of design variables in the RL approach at each iteration, in the \textsc{CTR} \& \textsc{DES} scenario. The RL algorithm converges at the final iteration to a battery capacity of approximately 44 kWh and a PV power output of around 57 kWp.}
    \label{fig:des_year_inv}
\end{figure}

From Table \ref{tab:results_year} it is also seen that the RL approach seems to promote solutions involving lower values of design variables. To further investigate the reasons underlying this result, the design variables for the evolution of the battery capacity and PV power, during the training process, are reported in Fig. \ref{fig:des_year_inv} in the \textsc{CTR} \& \textsc{DES} scenario.

As indicated in Table \ref{tab:parameters}, the design variable values can range from 0 to 200. However, it can be seen that higher values are not explored by the RL method. This latter resulted, at the last iteration, in design variables of 44 kWh for battery capacity and 57 kWp for PV power. During the training phase, the maximum values reached were 59 kWh and 90 kWp for battery capacity and PV power, respectively. This maximal explored battery capacity value is lower than the optimal one found by the MILP approach (95 kWh). Thus, the RL solution of the PV power value is expected to be lower. Indeed, the reward value is penalized if the RL agent injects PV production into the grid, since the cost of exported energy into the grid ($C_{grid}^{exp}$) is defined as a negative value in Table \ref{tab:parameters}. Consequently, limited battery capacity intrinsically causes a lower PV power value.

\section{Discussion}
\label{sec:discuss}

This section discusses the main observations that can be drawn from solving a battery-PV system with both RL and MILP approaches using the one-week and one-year data set.\\

\subsection{The promises of RL for joint design and control of energy systems}

The motivation for this study was to explore the potential of RL to enable joint control and design of energy systems. Tables \ref{tab:results_week} (one-week data) and \ref{tab:results_year} (one-year data) show that RL provides a solution that is close to the optimal MILP one. This is encouraging as it suggests that despite RL relying on a different optimization strategy, it is able to identify a meaningful solution in a simple case. However, the difference of reward value between MILP and RL increases when integrating design variables to the optimization problem, i.e., \textsc{CTR} \& \textsc{DES} scenarios in Tables \ref{tab:results_week} and \ref{tab:results_year}. Interestingly, the solutions for design variables are consistently smaller in RL as compared to MILP. Furthermore, from Fig. \ref{fig:des_year_inv}, it can be seen that the RL algorithm did not explore higher design variable values in the one-year case study. This observation can be explained by two possibilities: first, DEPS is a local-search method that is thus subject to converging towards local extrema. Once the control policy is too specialised to the investment parameters (under optimization too), these parameters are thus expected to be locally optimal and the algorithm is stuck. Second, the RL algorithms is subject to many hyperparameters to which the final results are sensible, it is possible that a different policy architecture, learning rate, or simply more iterations would ameliorate the performance of the method. Supporting the first explanation is the similarity between the reward values of the RL (-250) and MILP, based on same investments, (-247) approaches for the \textsc{CTR} \& \textsc{DES} scenario with $T=8760$ (Table \ref{tab:results_year}). Hence, in this specific energy system case study, it could be likely that the RL algorithm did not deem it advantageous to enhance the value of the design variables for either one or both of the two reasons stated.

Overall, these results show that RL provides realistic control and design strategies. Based on this, RL could be used to define new real-time control strategies integrating design constraints, and that are less sensitive to linearization inaccuracies \citep{reus_modeling_2019, sanchez_optimised_2019}. Given the differences in how uncertainties are accounted for by both methods, RL could also be a better candidate to integrate resources coming with high levels of uncertainty such as electric mobility.\\

\subsection{Technical challenges and future directions}

The main technical challenges encountered in this study are essentially the ones inherent to RL methods. First, various parameters need to be tuned: neural network architecture for the policy, the batch size for the optimization, the learning rate, or the different scaling among others. These parameters were tuned by trial and error and would need to be adapted to each new application. For example, the number of layers required in the one-year case study was larger than for the one-week toy example. Second, convergence of the RL method is not guaranteed, and when convergence happens, the solution is not guaranteed to be globally optimal. Third, as illustrated here above for the results of Figure \ref{fig:des_year_inv}, determining the number of iterations (set to 100,000 for the training phase in all our experiments) is also crucial and might affect RL solutions. Therefore, comparing RL and MILP solutions is not trivial because its is difficult to compare perfect foresight with policy based decisions. This should be accounted for when analyzing results from Tables \ref{tab:results_week} and \ref{tab:results_year}.\\

From a technical point of view, future work will aim at using more advanced RL methods. In particular, the RL algorithm used here is a modified version of the REINFORCE algorithm \citep{williams_simple_1992}, which was developed in 1992 and is one of the earliest RL algorithms. Today, more advanced algorithms are available for control problems, which can converge more rapidly or account for infinite time horizons, such as actor-critic algorithms (e.g., PPO \citep{schulman_proximal_2017} and GAE \citep{schulman_high-dimensional_2018}), but are yet to be adapted to joint design and control. In terms of applications, future work will aim to better evaluate the added value of RL by assessing the long-term performance of real-time sized systems. For example, a control framework could be developed to establish an operation strategy for the MILP-sized system. The framework would then be evaluated using several years of real-time data from the same system used for design. The same exercise would be applied to the trained model of the DEPS algorithm and performance obtained from several years of system control would be benchmarked, and the impact of design decisions could be discussed with more perspective.

\section{Conclusions}
\label{sec:ccl}

 In most studies, MILP is used for the design of energy systems and RL for the control. On the one hand, MILP assumes a perfect foresight of the future and is difficult to generalize to new data. On the other hand, RL methods proved to be efficient in other tasks linked to design and control but not on energy systems. In this study, we assessed the potential of an RL method, DEPS, i.e. an RL algorithm proven efficient for designing and controlling complex systems, for the joint design and control of energy systems. 

The energy system studied is a PV-battery system used to answer a real-life demand in order to minimize the overall cost. In order to assess the efficiency of the RL method, we compared the outcomes with those obtained with a MILP. As these two approaches are fundamentally different, the optimization problem was formulated in two distinct ways: first as a MILP and second as an MDP. The methodology and experimental context were clarified to facilitate the discussion of results and have a fair comparison. Both approaches are discussed in terms of their strengths and weaknesses. 

The findings show that RL can produce control strategies that are close to optimal, while using different values of design variables. This highlights the potential of RL for joint design and control of energy systems, particularly in scenarios where stochasticity is a key factor. However, the study also highlights the difficulty of tuning and using theses methods. Moving forward, there are several challenges to address, including the need to ensure that the RL solution converges to a global optimum. However, the promising results obtained in this study suggest that RL has the potential to be a valuable tool for jointly designing and controlling energy systems.


\newpage
\printbibliography

\end{document}